\documentclass[12pt]{article}
\usepackage{amsfonts,amsmath,amsxtra}
\usepackage{latexsym}
\usepackage{amssymb}
\usepackage{color}
\usepackage[matrix,arrow,curve]{xy}

\makeatletter
\def\@seccntformat#1{\csname
the#1\endcsname\enspace} \makeatother

\def\hybrid{\topmargin 0pt      \oddsidemargin 0pt
        \headheight 0pt \headsep 0pt
        \textwidth 16.5cm
        \textheight 23cm
        \marginparwidth 0.0in
        \parskip 5pt plus 1pt   \jot = 1.5ex}
\catcode`\@=11
\def\marginnote#1{}
\newcount\hour
\newcount\minute
\newtoks\amorpm
\hour=\time\divide\hour by60 \minute=\time{\multiply\hour by60
\global\advance\minute by-\hour}
\edef\standardtime{{\ifnum\hour<12 \global\amorpm={am}%
        \else\global\amorpm={pm}\advance\hour by-12 \fi
        \ifnum\hour=0 \hour=12 \fi
      \number\hour:\ifnum\minute<10 0\fi\number\minute\the\amorpm}}
\edef\militarytime{\number\hour:\ifnum\minute<10 0\fi\number\minute}

\def\draftlabel#1{{\@bsphack\if@filesw {\let\thepage\relax
   \xdef\@gtempa{\write\@auxout{\string
      \newlabel{#1}{{\@currentlabel}{\thepage}}}}}\@gtempa
   \if@nobreak \ifvmode\nobreak\fi\fi\fi\@esphack}
        \gdef\@eqnlabel{#1}}
\def\@eqnlabel{}
\def\@vacuum{}
\def\draftmarginnote#1{\marginpar{\raggedright\scriptsize\tt#1}}

\def\draft{\oddsidemargin -0.1truein
        \def\@oddfoot{\sl HeckeBaxter.tex \hfil
        \rm\thepage\hfil\sl\today\quad\militarytime}
        \let\@evenfoot\@oddfoot \overfullrule 3pt
        \let\label=\draftlabel
        \let\marginnote=\draftmarginnote
\def\@eqnnum{{\rm (\theequation)}
\rlap{\kern\marginparsep\tt\@eqnlabel}%
\global\let\@eqnlabel\@vacuum}  }

\newfont{\Bbbb}{msbm7 scaled 1\@ptsize00}
\newcommand{\zs}{\raise-1pt\hbox{$\mbox{\Bbbb Z}$}}

scaled 1\@ptsize00

\def\numberbysection{\@addtoreset{equation}{section}
        \def\theequation{\thesection.\arabic{equation}}}
\numberbysection

\renewcommand{\theequation}{\thesection.\arabic{equation}}

\def\titlepage{\@restonecolfalse\if@twocolumn\@restonecoltrue\onecolumn
     \else \newpage \fi \thispagestyle{empty}\c@page\z@
\def\thefootnote{\fnsymbol{footnote}} }
\def\endtitlepage{\if@restonecol\twocolumn \else  \fi
        \def\thefootnote{\arabic{footnote}}
        \setcounter{footnote}{0}}  
\relax \hybrid
\makeatletter
\newdimen\normalarrayskip            
\newdimen\minarrayskip               
\normalarrayskip\baselineskip \minarrayskip\jot
\newif\ifold             \oldtrue            \def\new{\oldfalse}
\def\arraymode{\ifold\relax\else\displaystyle\fi}
\def\eqnumphantom{\phantom{(\theequation)}} 
\def\@arrayskip{\ifold\baselineskip\z@\lineskip\z@
     \else
     \baselineskip\minarrayskip\lineskip1\baselineskip\fi}
\def\@arrayclassz{\ifcase \@lastchclass \@acolampacol \or
\@ampacol \or \or \or \@addamp \or
   \@acolampacol \or \@firstampfalse \@acol \fi
\edef\@preamble{\@preamble
  \ifcase \@chnum
     \hfil$\relax\arraymode\@sharp$\hfil
     \or $\relax\arraymode\@sharp$\hfil
     \or \hfil$\relax\arraymode\@sharp$\fi}}
\def\@array[#1]#2{\setbox\@arstrutbox=\hbox{\vrule
     height\arraystretch \ht\strutbox
     depth\arraystretch \dp\strutbox
width\z@}\@mkpream{#2}\edef\@preamble{\halign \noexpand\@halignto
\bgroup \tabskip\z@ \@arstrut \@preamble \tabskip\z@ \cr}%
\let\@startpbox\@@startpbox \let\@endpbox\@@endpbox
  \if #1t\vtop \else \if#1b\vbox \else \vcenter \fi\fi
  \bgroup \let\par\relax
  \let\@sharp##\let\protect\relax
  \@arrayskip\@preamble}
\def\eqnarray{\stepcounter{equation}%
              \let\@currentlabel=\theequation
              \global\@eqnswtrue
              \global\@eqcnt\z@
              \tabskip\@centering              
              \let\\=\@eqncr
              $$%
            \halign to \displaywidth  \bgroup
             \eqnumphantom \@eqnsel
      \hskip\@centering                               
    $\displaystyle  \tabskip\z@ {##}$%
    &\global\@eqcnt\@ne \hskip 2\arraycolsep
         $ \displaystyle  \arraymode{##}$\hfil
    &\global\@eqcnt\tw@ \hskip 2\arraycolsep
         $\displaystyle\tabskip\z@{##}$\hfil
         \tabskip\@centering
    &{##}\tabskip\z@\cr}
\makeatother

\def\IC{\mathbb{C}}

\def\IQ{\mathbb{Q}}
\def\IR{\mathbb{R}}
\def\IZ{\mathbb{Z}}

\def\CF {\mathcal{F}}

\def\CH {\mathcal{H}}

\def\CP {\mathcal{P}}

\def\CU {\mathcal{U}}
\def\CV {\mathcal{V}}
\def\CW {\mathcal{W}}

\def\g {{\gamma}}
\def\s {{\sigma}}

\def\e{\epsilon}

\def\pr {\partial}

\def\wt{\widetilde}
\def\wh{\widehat}

\def\End{{\rm End}}

\def\Ind{{\rm Ind}}
\def\Irr{{\rm Irr}}
\def\Fun{{\rm Fun}}

\def\Lie{{\rm Lie}}
\def\Mat{{\rm Mat}}

\def\RE{{\mathop{\rm Re}}}

\def\sign{{\mathop{\rm sign}}}
\def\Span{{\mathop{\rm span}}}

\def\frak{\mathfrak}
\def\ov {{\overline}}

\def\Tr{{\rm Tr}\,}

\def\<{\langle}
\def\>{\rangle}
\def\ov{\overline}

\makeatletter
\DeclareRobustCommand{\loplus}{\mathbin{\mathpalette\dog@lsemi{+}}}
\DeclareRobustCommand{\roplus}{\mathbin{\mathpalette\dog@rsemi{+}}}

\newtheorem{te}{Theorem}[section]
\newtheorem{de}{Definition}[section]
\newtheorem{prop}{Proposition}[section]           
\newtheorem{cor}{Corollary}[section]
\newtheorem{lem}{Lemma}[section]

\newcommand{\proof}{\noindent {\it Proof}. }
\newcommand\bqa{\begin{eqnarray}}
\newcommand\eqa{\end{eqnarray}}
\def\be{\begin{eqnarray}\new\begin{array}{cc}}
\def\ee{\end{array}\end{eqnarray}}
\def\beq{\begin{equation}}
\def\eeq{\end{equation}}
\def\bse{\begin{subequations}}                
\def\ese{\end{subequations}}
\def\bp{\begin{pmatrix}}
\def\ep{\end{pmatrix}}

\def\i{\imath}
\newcounter{pac}[section]

\newcounter{pacc}[subsection]

 0

\begin{document}

\title{\bf The $GL_{\ell+1}(\IR)$ Hecke-Baxter operator:\\
principal series representations}
\author{A.A. Gerasimov, D.R. Lebedev and S.V. Oblezin}
\date{\today}
\maketitle

\renewcommand{\abstractname}{}

\begin{abstract}
 \noindent {\bf Abstract}.
Previously introduced the $GL_{\ell+1}(\IR)$ Hecke-Baxter operator
is a one-parameter family of elements in the commutative spherical
Hecke algebra\\ $\CH(GL_{\ell+1}(\IR),O_{\ell+1})$. Its action  on
spherical vectors in  spherical principle series representations of
$GL_{\ell+1}(\IR)$ is given by multiplication by the Archimedean
$L$-factors associated to these representations. In this note we
propose  an extension of the  construction to other (non-spherical)
$GL_{\ell+1}(\IR)$ principle series representations providing a
relevant generalization of the notions of spherical vector,
commutative spherical Hecke algebra and the Hecke-Baxter operator to
the general case. Action of the introduced Hecke-Baxter operator on
the generalized spherical vectors  is given by multiplication by the
Archiemdean $L$-factor associated to the corresponding principle
series representation of $GL_{\ell+1}(\IR)$.
\end{abstract}

 \vspace{5mm}


\section{Introduction}

The local Langlands correspondence associates to representations of
the Galois group of a local field $K$ (or its proper substitute)
admissible representations of reductive Lie groups over $K$. In more
concrete terms one associates to an admissible $G(K)$-representation
a function in one complex variable $s$ called local $L$-factor.
These functions capture invariants of the corresponding
representation of the Galois group of $K$. The simplest instance of
local $L$-factors  appears in the product decomposition   of the
(completed) Riemann zeta-function associated in this setup to
representation theory of the group $GL_1$ over the local completions
of $\IQ$. Starting with the Tate construction \cite{T} for $GL_1$
and further generalization to the case of $GL_2$ in \cite{JL}, the
corresponding $L$-factors are constructed in a rather indirect way
by invoking Fourier transformations of certain class of measures on
reductive Lie groups over local fields. Obviously, a more direct
approach would be desirable.

In \cite{GLO08} we have considered this problem in the special case
of spherical principal series representations of $GL_{\ell+1}(\IR)$.
We have proposed a one-parameter family of elements of the
commutative spherical Hecke algebra
$\CH(GL_{\ell+1}(\IR),O_{\ell+1})$ such that its action on uniquely
defined spherical vectors in  spherical principal series
representations of $GL_{\ell+1}(\IR)$ is given by multiplication by
the corresponding local Archimedean $L$-factor. This construction
was motivated by considerations in the theory of quantum integrable
systems, specifically the $GL_{\ell+1}(\IR)$-Toda chain \cite{GLO08}
(see also \cite{G}).
Note that this construction naturally incorporates the fact that
local Archimedean $L$-factors have structure of $\IR^*_+$-torsors
(to fix a section of this torsor one should invoke global arithmetic
considerations).

In this note we propose an extension of the construction of the
Hecke-Baxter operator to the case of general (non-spherical)
principal series representations  of $GL_{\ell+1}(\IR)$. We start
with introducing special vectors in a principal series
representations generalizing  in a sense  the spherical vector to
the case of non-spherical principal series representations. Next, we
define a still commutative extension
$\CH_r(GL_{\ell+1}(\IR),O_{\ell+1})$  of the standard spherical
Hecke algebra $\CH(GL_{\ell+1}(\IR),O_{\ell+1})$.  Finally,  we
introduce a one-parameter family of elements in
$\CH_r(GL_{\ell+1}(\IR),O_{\ell+1})$ considered as a generalization
of the Hecke-Baxter operator proposed in \cite{GLO08}. By
construction, the Hecke-Baxter operators for various values of the
parameter commute and thus one can look for common eigenvectors of
these operators in principal series representations. The main result
of this note is given by Theorem \ref{MAINTHM} claiming  that the
action of the proposed Hecke-Baxter operator on the generalized
spherical vector in a principal series representation reduces to
multiplication by the corresponding local Archimedean $L$-factor.
Using the notion of the generalized spherical vector we  introduce
the corresponding generalized  spherical and Whittaker functions as
certain matrix elements in  principal series representations. The
proven Theorem \ref{MAINTHM} entails that the action of the
Hecke-Baxter operator on these matrix elements is given by
multiplication the $L$-factors attached to the representations of
$GL_{\ell+1}(\IR)$.

Consideration of the non-spherical principal series reveals an
interesting subtlety in the definition of local Archimedean
$L$-factors. The standard integral representation of $L$-factors is
given in terms of integrals over Gaussian type measures on
non-compact Lie groups. The fact that these measures are invariant
under  Fourier transform plays a key role in deriving functional
equations  for global zeta-functions constructed from local
contributions given by local $L$-factors. Thus it is natural to look
for a compatibility of the introduced Hecke-Baxter operator acting
uniformly in various principal series representations of
$GL_{\ell+1}(\IR)$  with the  Fourier transform. It turns out that
the transformation properties are most simple  if the corresponding
Gaussian measure is replaced by the Feynmann measure (or quadratic
character in terms of \cite{Wei}). As usual, the latter  is
considered as a limit of a certain complex measure of the Gaussian
type. This opens an interesting possibility to look at number
theoretic zeta-functions (and their local counterparts like local
$L$-factors) in terms of  a kind of Quantum Field Theory over
arithmetic varieties  \cite{Ma}. This is obviously an interesting
direction to pursue. A less ambitious direction (although, not quite
unrelated) would be to provide an interpretation of the extended
Hecke-Baxter operator in terms two-dimensional topological quantum
field theories on the disk (i.e. in terms of brane geometry) along
the line of \cite{GLO11}. Finally let us stress an obvious
importance of the generalization of the proposed constructions  to
other series of admissible representations of $GL_{\ell+1}(\IR)$ as
well as its counterparts for $GL_{\ell+1}(\IC)$  that we are going
to report elsewhere.

{\it After submission of  this paper to arXiv we were informed by
  P. Humphries about similar results concerning generalization of
  the Hecke-Baxter formalism \cite{GLO08}
 to non-spherical   principal series representations of $GL_{\ell+1}(\IR)$ (see
  \cite{Hum} and references therein).
The construction of \cite{Hum} relies on an  extension of the
recursive relations for the spherical $GL_{\ell+1}(\IR)$-Whittaker
functions \cite{GLO08} to the non-spherical case \cite{Lin}. Yet
various results and approaches of \cite{Hum} and of the present
paper are complimentary and thus worth of further comparison.}

{\it  Acknowledgements:} The research of S.V.O. is partially
supported by the Beijing Natural Science Foundation grant IS24004.

\section{Spherical $GL_{\ell+1}(\IR)$ Hecke-Baxter operator}

In this Section we recall the known constructions of the spherical
Hecke algebras and the Hecke-Baxter operators for the spherical principal
series representations of the Lie group $GL_{\ell+1}(\IR)$
\cite{GLO08}. This allows us to proceed in the following Section with a
generalization of these constructions to the case of non-spherical
(ramified) principal series representations.

In \cite{GLO08}, the Hecke-Baxter operator was defined as a
one-parameter family of elements of the spherical Hecke algebra
$\CH(GL_{\ell+1},O_{\ell+1})$ associated with the Gelfand pair
$(GL_{\ell+1}(\IR),O_{\ell+1})$ where $O_{\ell+1}$ is the orthogonal
subgroup of $GL_{\ell+1}(\IR)$. The spherical Hecke algebra is a
commutative associative algebra modeled on the appropriate subspace
of $O_{\ell+1}$-biinvariant functions in
$L^1\bigl(GL_{\ell+1}(\IR)\bigr)$ with the algebra structure defined
by the convolution
 \be\label{conv}
  (f_1*f_2)(\tilde{g})\,=\!\int\limits_{GL_{\ell+1}(\IR)}\!\!
  d\mu(g)\,\,f_1(g)\,\,f_2(g^{-1}\tilde{g})\,.
 \ee
 Here $\mu(g)$ is the Haar measure on $GL_{\ell+1}(\IR)$
 \be\label{Hmes}
  d\mu(g)\,=\,
  |\det g|^{-(\ell+1)}\prod_{i,j=1}^{\ell+1}dg_{ij}\,,
 \ee
and conditions on the class of considered functions are imposed to
render the convolution operation to be defined.

Let $B\subset GL_{\ell+1}(\IR)$ be  the Borel subgroup identified
with the subgroup of lower-triangular matrices. Define the following
character of $B$,
 \be\label{Bchar}
  \chi^B_{\epsilon,\gamma}(b)\,
  =\prod_{j=1}^{\ell+1}\sign(b_{jj})^{\epsilon_j}\,
  |b_{jj}|^{\imath\gamma_j+\rho_j},\quad
  \gamma_j\in \IR,\quad \epsilon_j\in \{0,1\},\quad
  \rho_j=\frac{\ell}{2}+1-j\,,
 \ee
 trivial on the unipotent radical $N\subset B$. Principal
series $GL_{\ell+1}(\IR)$-representation
$(\pi_{\epsilon,\gamma},\CV_{\epsilon,\gamma})$ is the induced
representation
$\CV_{\epsilon,\gamma}=\Ind_B^{GL_{\ell+1}(\IR)}\chi^B_{\epsilon,\gamma}$
realized in the space of $B$-equivariant functions satisfying
 \be\label{PSrep}
  \phi(gb)=\chi^B_{\epsilon,\gamma}(b)\,\phi(g)\,,\qquad b\in B\,.
 \ee
The group $GL_{\ell+1}(\IR)$ acts on $\phi\in\CV_{\epsilon,\gamma}$
via the left action:
 \be
  \bigl(\pi_{\epsilon,\gamma}(g)\,\phi\bigr)(\tilde{g})\,
  =\,\phi(g^{-1}\tilde{g})\,,\quad g,\tilde{g}\in
  GL_{\ell+1}(\IR)\,.
 \ee
The representation
$(\pi_{\epsilon,\gamma},\CV_{\epsilon,\gamma})$ is irreducible for a
generic $\gamma$. The space
$\CV_{\epsilon,\gamma}=\Ind_B^{GL_{\ell+1}(\IR)}\chi^B_{\epsilon,\gamma}$
can be supplied with an invariant Hermitian form
$\<\,,\,\>$.

The principal series representation
$(\pi_{\epsilon,\gamma},\CV_{\epsilon,\gamma})$ is called spherical
if  the corresponding character $\chi^B_{\epsilon,\gamma}$ is
trivial on the subgroup $M=O_{\ell+1}\cap B$, or equivalently, if
$\e_j=0,\,1\leq j\leq\ell+1$ in \eqref{Bchar}. Each spherical
representation  contains a unique spherical vector $\phi_0$
invariant under the action of the subgroup $O_{\ell+1}$ and
normalized by the condition $\phi_0(1)=1$ (see e.g \cite{He}).

There are deep arithmetical reasons  to associate to each principal
series representation $\CV_{\epsilon,\gamma}$ a function in an
auxiliary variable $s\in \IC$ called the Archimedean $L$-factor (see
e.g. \cite{Kna} and Appendix in \cite{Ja}):
 \be\label{Lf0}
  L(s|\epsilon,\gamma)\,
  =\prod_{j=1}^{\ell+1}
  \pi^{-\frac{s+\epsilon_j-\imath\gamma_j}{2}}\,
  \Gamma\left(\frac{s+\epsilon_j-\imath\gamma_j}{2}\right)\,.
  \ee
Notice that the  $L$-factors  associated to principal series
$GL_{\ell+1}(\IR)$-representations for $\epsilon\neq 0$  enter product decomposition
of the $L$-factors associated with  spherical  principal series representations of the
complex group  $GL_{\ell+1}(\IC)$. For example  the following
  product  relation holds
(see e.g. Appendix in \cite{Ja}):
 \be
  L^{GL_{\ell+1}(\IC)}(s|\gamma)= L(s|0,\gamma) \cdot
  L(s|\epsilon_*,\gamma)\,,\qquad \epsilon_*=(1,\ldots, 1)\,,
 \ee
generalizing the standard Legendre identity
 \be
  \Gamma(s)\,
  =2^{s-1}\,\pi^{-\frac{1}{2}}\,\Gamma\Big(\frac{s}{2}\Big)\,
  \Gamma\Big(\frac{s+1}{2}\Big)\,.
  \ee
It is useful to consider more general local Archimedean
  $L$-factors  depending on $c\in \IR_+^*$
  \be\label{Lf}
  L(s,c|\epsilon,\gamma)\,
  =\prod_{j=1}^{\ell+1}
  c^{-\frac{s+\epsilon_j-\imath\gamma_j}{2}}\,
  \Gamma\left(\frac{s+\epsilon_j-\imath\gamma_j}{2}\right)\,.
  \ee
The choice of $c=\pi$ in \eqref{Lf0} is dictated by global
arithmetic considerations which are not relevant (except the final
Section 7) to our local
considerations of the representation theory of Lie groups defined
over $\IR$. Thus in the following we
consider \eqref{Lf} with some fixed non-specified $c\in\IR^*_+$\,, and
we use the simplified notation $L(s|\epsilon,\gamma)$.

For the spherical principal series representation
 $(\pi_{0,\gamma},\CV_{0,\gamma})$, consider
representation of the convolution algebra
$\bigl(L^1\bigl(GL_{\ell+1}(\IR)\bigr),\,*\bigr)$ in $\CV_{0,\gamma}$
defined by the action via the translation
operators. Namely, for $f\in L^1\bigl(GL_{\ell+1}(\IR)\bigr)$ we take
 \be\label{GArep}
  f\cdot \phi\,:
  =\!\int\limits_{GL_{\ell+1}(\IR)}\!\!
  d\mu(g)\,\,f(g)\,\,\,\pi_{0,\gamma}(g)\cdot
  \phi\,,\qquad\phi\in\CV_{0,\gamma}\,,
  \ee
provided the corresponding integral is convergent. In particular the
action is defined for  elements of the spherical Hecke algebra
$\CH(GL_{\ell+1}(\IR),O_{\ell+1})$ represented by
$O_{\ell+1}$-biinvariant functions on $GL_{\ell+1}(\IR)$. Uniqueness
of the spherical vector $\phi_0$ in the spherical representation
$(\pi_{0,\gamma},\CV_{0,\gamma})$ gives rise to the action of
elements of  $\CH(GL_{\ell+1}(\IR),O_{\ell+1})$ on $\phi_0$ via  one-dimensional
representation $\Lambda$ of the commutative algebra
$\CH(GL_{\ell+1}(\IR),O_{\ell+1})$
 \be\label{Eig23}
  f\cdot \phi_0\,
  =\,\Lambda_f\,\phi_0\,,\quad f\in
  \CH(GL_{\ell+1}(\IR),O_{\ell+1})\,,\quad\Lambda_f\in\IC\,.
 \ee
Thus the spherical
vector $\phi_0\in \CV_{0,\gamma}$ is a common
$\CH(GL_{\ell+1}(\IR),O_{\ell+1})$-eigenvector.

The Iwasawa decomposition  of the Lie group $GL_{\ell+1}(\IR)$ has the
following form
 \be\label{Iwasawa}
  GL_{\ell+1}(\IR)\,=\,O_{\ell+1}\,A N\,,
  \ee
  where  for our choice of the Borel subgroup $B$, $A$
  is the group of the diagonal matrices with strictly positive
 real entries. Let us introduce the corresponding product representation for element
 of $GL_{\ell+1}(\IR)$
 \be\label{Iww}
 g=kan, \qquad g\in GL_{\ell+1}(\IR), \qquad
 k\in O_{\ell+1},\qquad a\in A,\qquad n\in N\,.
 \ee
 Using  \eqref{Bchar} the spherical vector $\phi_0\in \CV_{0,\gamma}$
might be written explicitly as
 \be\label{SphVect}
  \phi_0(kan)\,
  =\,\chi^B_{\gamma}(an)\,
  =\prod_{j=1}^{\ell+1}a_j^{\imath\gamma_j+\rho_j}\,.
  \ee
In the following we will use a more concise notation
 $\chi^B_{\gamma}(b):= \chi^B_{0,\gamma}(b)$.

 In \cite{GLO08}, we introduced a one-parameter family of elements
 \be
 Q_{s}\in\CH(GL_{\ell+1}(\IR),O_{\ell+1}),\,\qquad s\in\IC\,,
\ee
such that its action on the spherical vector $\phi_0$ in a spherical principal series
$GL_{\ell+1}(\IR)$-representation $(\pi_{0,\gamma},\CV_{0,\gamma})$
reduces to multiplication  by  the corresponding $L$-factor
\eqref{Lf} with $\epsilon=0$. Such  $Q_s$ was  called the Hecke-Baxter
operator associated to $(\pi_{0,\gamma},\CV_{0,\gamma})$.

\begin{prop}\label{SphHB}
The spherical vector $\phi_0\in\CV_{0,\gamma}$
  is the eigenfunction of the
  Hecke-Baxter operator acting on  by convolution with the following
function on $GL_{\ell+1}(\IR)$
 \be\label{HB0}
  Q_{s}(g)\,
  =\,(c\pi^{-1})^{\frac{\ell(\ell+1)}{4}}\,
  |\det g|^{s+\frac{\ell}{2}}\,\,
  e^{-c\Tr(g^{\top}g)}\,,\qquad c\in\IR^*_+\,,\qquad {\rm Re}(s)>0\,.
 \ee
The corresponding eigenvalue is given by the Archiemdean $L$-factor attached to
the representation $(\pi_{0,\gamma},\CV_{0,\gamma})$:
 \be\label{Lf10}
  L(s|0,\gamma)\,
  =\prod_{j=1}^{\ell+1}c^{-\frac{s-\imath\gamma_j}{2}}\,
  \Gamma\left(\frac{s-\imath \gamma_j}{2}\right)\,.
 \ee
\end{prop}

\proof The function $Q_{s}(g)$ is an element of
$\CH(GL_{\ell+1}(\IR),O_{\ell+1})$, and thus  according to
\eqref{Eig23}  the spherical vector \eqref{SphVect} is its
eigenfunction with an eigenvalue  $\Lambda(s|\gamma)$:
 \be
  (Q_{s}*\phi_0)(\tilde{g})\,
  =\int\limits_{GL_{\ell+1}(\IR)}\!\!d\mu(g)\,\,
  Q_{s}(g)\,\,\phi_{0}(g^{-1}\tilde{g})\,
  =\,\Lambda(s|\gamma)\,\,\phi_0(\tilde{g})\,.
 \ee
Taking into account the normalization condition $\phi_{0}(1)=1$ for
the spherical vector \eqref{SphVect} implies the following
integral expression for the eigenvalue
 \be
  \Lambda(s|\gamma)\,
   =\!\int\limits_{GL_{\ell+1}(\IR)}\!\!d\mu(g)\,\,
  Q_{s}(g^{-1})\,\,\phi_{0}(g)\,,
  \ee
where the unimodularity of $GL_{\ell+1}(\IR)$ is used.

Under  the Iwasawa decomposition \eqref{Iwasawa}, the  Haar measure
\eqref{Hmes} on $GL_{\ell+1}(\IR)$  allows the decomposition as a
product of the Haar measures on $A$, $O_{\ell+1}$ and $N$
 \be\label{Haardec}
  d\mu(g)=\delta_N(a)\,dk\times da\times dn\,,\qquad g=kan\,,
 \ee
where  $\delta_N(a)=a^{-2\rho}$ is the modular function of the
unipotent subgroup $N\subset B$. The corresponding Haar measures on the subgroups are
given by
 \be
  dn=\prod_{i>j}^{\ell+1} dn_{ij},\qquad
  da=\prod_{i=1}^{\ell+1} \frac{da_i}{a_i}\,,
 \ee
and the Haar measure
$dk$ on $O_{\ell+1}$ is normalized by the condition \be\label{Knorm}
\int_{O_{\ell+1}}dk=1\,. \ee Here we include  the sum over connected
components of $O_{\ell+1}$ into the integral. Applying
\eqref{Haardec}  we derive
 \be
  \Lambda(s|\gamma)\,
  =\!\int\limits_{O_{\ell+1}}\! dk
  \int\limits_{A}\!da\int\limits_N\!dn\,\,\delta_{N}(a)\,\phi_{0}(a)\,\,
  Q_{s}(n^{-1}a^{-1})\,.
 \ee
 Substituting $\eqref{HB0}$ and
taking into account \eqref{Knorm} we
 arrive at  the following integral:
 \be\label{EIGENVAL}
  \Lambda(s|\gamma)\,
  =\,(\pi c^{-1})^{-\frac{\ell(\ell+1)}{4}}\!
  \int\limits_{A}\!da\,\,
  a^{-2\rho}\,\chi^B_{\g}(a)\,\,
  |\det a|^{-s-\frac{\ell}{2}}\\
  \times\int\limits_N\!dn\,\,e^{-c\Tr(n^{\top}a^{-2}n)}\,.
 \ee
Next, we note that
 \be
  \Tr\bigl(n^\top a^{-2}n\bigr)\,
  =\,\sum_{i=1}^{\ell+1}a_i^{-2}\,
  +\,\sum_{i>j}n^2_{ij}a_i^{-2}\,.
 \ee
Then using the standard Gauss integral formula for the integration
over $n\in N$ gives
 \be
  \Lambda(s|\gamma)\,
  =\!\int\limits_{A}\!da\,\,
  |\det a|^{-s-\frac{\ell}{2}}
  \prod_{j=1}^{\ell+1}a_j^{\i\g_j-\rho_j+\ell+1-j}\,
  e^{-\frac{c}{a_i^2}}\\
  =\prod_{j=1}^{\ell+1}
  \int\limits_0^{\infty}\frac{da_j}{a_j}\,\,
  a_j^{s+\frac{\ell}{2}-\i\g_j+\rho_i-\ell-1+j}\,\,
  e^{-ca_j^2}\,
  =\prod_{j=1}^{\ell+1}c^{-\frac{s-\imath\gamma_j}{2}}\,
  \Gamma\left(\frac{s-\imath \gamma_j}{2}\right)\,,
  \ee
where in the latter equality we take into account the condition
$\RE(s)>0$ to apply the Euler
integral formula for the Gamma-function. $\Box$

Now, let $N_+\subset GL_{\ell+1}(\IR)$ be the maximal unipotent
subgroup opposite to $N\subset B$ and let $\chi^{N_+}:N_+\to\IC^*$
be a principal character. Define the Whittaker vector
$\psi\in\CV_{0,\gamma}$ by the following condition:
 \be
  \pi_{0,\gamma}(n)\cdot \psi=\chi^{N_+}(n)\,\psi\,,\qquad n\in N_+\,.
  \ee
By \cite{Sha}, the Whittaker vector is unique in $\CV_{0,\gamma}$.
The $GL_{\ell+1}(\IR)$-invariant Hermitian pairing $\<\,,\,\>$
on the spherical principal series representation
$(\pi_{0,\g},\CV_{0,\gamma})$ allows  to define  the spherical and Whittaker
  functions   associated to the
principal series representation $(\pi_{0,\g},\CV_{0,\gamma})$ as the
following matrix elements:
 \be\label{SPFUN}
  \Phi_{0,\gamma}(g)=\<\phi_0,\,\pi_{0,\gamma}(g)\,\phi_{0}\>\,, \qquad
  \Psi_{0,\gamma}(g)=\<\phi_0,\,\pi_{0,\gamma}(g)\,\psi\>\,.
 \ee

\begin{cor}
The Hecke-Baxter operator acting on the spherical and Whittaker
functions \eqref{SPFUN} by convolution with the function
 \be
  Q_{s}(g)\,
  =\,(c \pi^{-1})^{\frac{\ell(\ell+1)}{4}}\,
  |\det g|^{s+\frac{\ell}{2}}\,\,
  e^{-c\Tr(g^{\top}g)}\,,
 \ee
has an eigenvalue given by the Archiemdean $L$-factor attached to
the representation
 \be
  L(s|0,\gamma)=\prod_{j=1}^{\ell+1}c^{-\frac{s-\imath\gamma_j}{2}}\,
  \Gamma\left(\frac{s-\imath \gamma_j}{2}\right)\,.
 \ee
\end{cor}

\section{Generalized spherical vectors}

In this Section we propose a  generalization of the notion of
the spherical vector to the case of general principal series representations
$(\pi_{\epsilon,\g},\CV_{\epsilon,\g})$ with $\epsilon\neq 0$.
 Let us start with considering the finite group
$M=O_{\ell+1}\cap B$ isomorphic to $\IZ_2^{\ell+1}$.  Explicitly
the group $M$ may be identified with  the   subgroup of diagonal matrices
with entries in $\{\pm 1\}$
 \be\label{Msub}
  m={\rm diag}\bigl((-1)^{\alpha_1}, \cdots ,
  (-1)^{\alpha_{\ell+1}}\bigr), \qquad \alpha_i\in \{0,1\}\,.
 \ee
Irreducible representations  of $M$ are  one-dimensional  and
 may be parameterized  by the set $\CP$ of signatures
 \be\label{SIGNA}
  \epsilon=(\epsilon_1, \cdots ,\epsilon_{\ell+1})\,,\qquad   \epsilon_i\in
  \{0,1\}\,.
  \ee
For each $\e\in\CP$, the corresponding  one-dimensional
$M$-representation $\chi^M_{\e}$ is given by
 \be\label{Mchar}
  \chi^M_\epsilon(m)=(-1)^{\sum\limits_{i=1}^{\ell+1}\alpha_i\epsilon_i},
 \ee
for $m\in M$ represented in the form \eqref{Msub}.
For a signature
$\e\in\CP$, define the module $|\epsilon|$ by
 \be
  |\epsilon|=\sum_{i=1}^{\ell+1}\epsilon_i\,\,.
 \ee

Let $V=\IC^{\ell+1}$ be the standard irreducible representation of
$GL_{\ell+1}(\IR)$ and let $\{e_i,\,1\leq i\leq\ell+1\}\subset V$ be
an orthonormal basis, $(e_i,\,e_j)=\delta_{ij}$. Let $(\pi_k,\,W_k)$
be the fundamental representations of $GL_{\ell+1}(\IR)$ in the spaces
 \be\label{Fundrep}
  W_k\,=\,\wedge^kV=\wedge^k\IC^{\ell+1},\quad
  k=0,1,\ldots,\ell+1\,,
 \ee
so in particular, $W_1=V$, and $W_0=\wedge^0V\simeq\IC$ is the
trivial representation.  For each $k=1,\ldots,\ell+1$, choose an
  orthonormal basis of poly-vectors in $W_k$
  parameterized by the signatures $\e$ with $|\epsilon|=k$
 \be\label{Wbasis}
  v_\epsilon\,
  =\,(k!)^{1/2}\cdot
e_1^{\epsilon_1}\wedge \cdots \wedge
  e_{\ell+1}^{\epsilon_{\ell+1}}\,,\qquad
  (v_{\e},\,v_{\e'})=\delta_{\e,\e'}\,,\qquad |\epsilon|=|\epsilon'|=k\,,
  \ee
where we use the following convention
 \be
  v_1\wedge \cdots \wedge v_n\,:
  =\,\frac{1}{n!}\sum_{\s\in \mathfrak{S}_n}
  v_{\s(1)}\otimes\cdots\otimes v_{\s(n)}\,.
 \ee
Then in the basis \eqref{Wbasis}, the action of $M\subset
O_{\ell+1}$ in $W_k$ has a simple form:
 \be\label{Mact}
  \pi_k(m)\,
  v_\epsilon\,
  =\,\chi^M_{\epsilon}(m)\,v_{\epsilon}\,,\qquad m\in M\,.
 \ee
The  basis \eqref{Wbasis} is uniquely determined as the eigenbasis
in each $W_k$ diagonalizing  the $M$-action. Combining
\eqref{Wbasis} for all $0\leq k\leq\ell+1$ results in the basis
$\{v_{\e},\,\e\in\CP\}$ of the (reducible)
$GL_{\ell+1}(\IR)$-representation $(\pi_W,W)$,
 \be\label{W}
  W\,=\,\bigoplus_{k=0}^{\ell+1}W_k\,,
 \ee
with the $M$-action given by \eqref{Mact}. We  may identify $W$
with the  regular representation of $M$. In the following we will
use the diagonal matrix elements in representation $W$. These matrix
elements might  be  expressed in terms of the principal minors
$\Delta_{\epsilon}(g)$ of the matrix $g$ with the rows and columns
corresponding to non-zero $\epsilon_i$'s
\be
  (v_{\epsilon},\,\pi_{W}(g)\,v_{\epsilon})\,
  =\,\Delta_{\epsilon}(g)\,,\qquad g\in GL_{\ell+1}(\IR)\,.
 \ee

The next simple statement will be instrumental in the following
constructions.

\begin{lem} Upon restriction
  to $O_{\ell+1}$, the fundamental representations
 $(\pi_k,\,W_k)$ of the group $GL_{\ell+1}(\IR)$ remain irreducible.
\end{lem}
\proof Recall that finite-dimensional irreducible representations of
$GL_{\ell+1}(\IR)$ are labeled by the Young diagrams and may be
realized in  tensor powers of the standard representation
$V=\IC^{\ell+1}$  \cite{Wey}. In particular, the fundamental
representations $(\pi_k,W_k)$ correspond to the $k$-column Young
diagrams and are given by  the totally skew-symmetric tensor
representations.  In turn, finite-dimensional irreducible
representations of $O_{\ell+1}$ are described by the same types of
tensors  with the additional condition on the tensors
to be traceless. The assertion follows by the fact that the
traceless condition is vacuous for the fundamental representations
of $GL_{\ell+1}(\IR)$.  $\Box$

Now  consider the  principal series representation
$(\pi_{\epsilon,\gamma},\CV_{\epsilon,\gamma})$ induced from the
 character of  the Borel subgroup $B \subset GL_{\ell+1}(\IR)$
 accordingly to \eqref{Bchar}. We would like to introduce
a vector $\phi_\epsilon\subset \CV_{\epsilon,\gamma}$ generalizing
in a sense the notion of the spherical vector $\phi_0\in \CV_{0,\g}$
in \eqref{SphVect}. To do
this we invoke spherical model for the principal series
representation $\CV_{\epsilon,\gamma}$. Precisely, using the Iwasawa
decomposition \eqref{Iwasawa} and
the presentation $B=MAN$ of the Borel subgroup  we infer from the
\eqref{PSrep} the identification of $\CV_{\epsilon,\gamma}$ with the space of
$M$-equivariant functions on $O_{\ell+1}$:
 \be\label{ModelK}
  \CV_{\epsilon, \gamma}\,\simeq\,\bigl\{\phi\in L^2(O_{\ell+1})\,:\quad
  \phi(km)=\chi^M_{\epsilon}(m)\,\phi(k),\,\, m\in M\bigr\}\,,
 \ee
with respect to the characters \eqref{Mchar},
 \be
  \chi^{M}_\epsilon(m)\,
  =\prod_{j=1}^{\ell+1}\sign(m_{jj})^{\epsilon_j}\,.
 \ee
Therefore, vectors in $\CV_{\e,\gamma}$ can be described explicitly using
the Peter-Weyl theory for compact groups. Indeed, for the Lie group
$O_{\ell+1}$, there is a canonical decomposition (of the left/right regular
representation) in terms of matrix elements of its irreducible
representations:
 \be\label{Kdec}
  L^2(O_{\ell+1})=\bigoplus_{\mu\in\Irr(O_{\ell+1})}\,
  (V^*_\mu \otimes V_\mu)\,.
 \ee
Then by \eqref{ModelK} the principal series representation
$\CV_{\epsilon,\gamma}$ allows the decomposition
 \be\label{Kdec1}
  \CV_{\epsilon,\gamma}
  =\bigoplus_{\mu\in\Irr(O_{\ell+1})}
  (V^*_\mu \otimes V_{\mu,\epsilon})\,,
 \ee
where $V_{\mu,\epsilon}\subset V_\mu$ is the subspace of vectors
such that the group $M$ acts on these vectors via character
$\chi^M_{\epsilon}$. Let $\CV^W_{\epsilon,\gamma}$
be the subspace of $\CV_{\e,\g}$ (realized  via  \eqref{Kdec1})
spanned by  matrix
elements of the irreducible representations of $O_{\ell+1}$ obtained
by restriction of the $GL_{\ell+1}$-representation $(\pi_W,W)$ in
\eqref{W}. Then the following matrix elements provide a basis in
$\CV^W_{\epsilon,\gamma}$ enumerated by $\e'\in\CP$:
 \be\label{Wsub}
  \phi^{\epsilon'}_{\epsilon,\gamma}(g)=
  (v_{\epsilon'},\,\pi_W(k)\,v_{\epsilon})\,
  \,\chi^B_{\gamma}(a)\,, \qquad g=kan\,,
  \ee
  where  the Iwasawa decomposition  of $g\in  GL_{\ell+1}(\IR)$ defined
  by \eqref{Iww} is used.
 Note that the above matrix elements are
non-trivial only if $|\epsilon'|=|\epsilon|$, and are normalized
accordingly to \eqref{Wbasis}:
 \be
  \phi^{\epsilon'}_{\epsilon,\gamma}(1)= (v_{\epsilon'},\,\pi_W(1)\,v_{\epsilon})\,
  =\,(v_{\e'},\,v_{\e})\,
  =\,\delta_{\epsilon',\epsilon}\,,\qquad\e,\e'\in\CP\,.
 \ee
The basis elements \eqref{Wsub} possess the following two obvious properties:
 \be\label{Wmatel}
  \phi^{\epsilon_1}_{\epsilon,\gamma}(kg)\,
  =\,\sum_{\epsilon_2\in\CP\atop|\e_2|=|\e_1|}\,
  (v_{\e_1},\,\pi_{W}(k)\,v_{\e_2})\,
  \phi^{\epsilon_2}_{\epsilon,\gamma}(g),
  \quad k\in O_{\ell+1}\,,\\
  \phi^{\epsilon_1}_{\epsilon,\gamma}(m_1gm_2)\,
  =\,\chi^M_{\epsilon_1}(m_1^{-1})\,\chi^M_{\epsilon}(m_2)\,
  \phi^{\epsilon_1}_{\epsilon,\gamma}(g)\,,\quad m_1,m_2\in M\,.
 \ee

\begin{de}\label{DEF}
Let $(\pi_{\epsilon,\gamma},\CV_{\epsilon,\gamma})$ be
  a principal series representation of $GL_{\ell+1}(\IR)$
induced from the  character of the Borel subgroup $B=MAN$,
 \be
  \chi^{B}_{\e,\gamma}(man)=\prod_{j=1}^{\ell+1}
  \sign(m_j)^{\epsilon_j}\,\,a_j^{\imath \gamma_j+\rho_j}, \qquad
  m\in M,\quad
  a\in A,\quad n\in N\,.
  \ee
Let $\CV^W_{\epsilon,\gamma}$ be subspace  of
$\CV_{\epsilon,\gamma}$  spanned by the matrix elements \eqref{Wsub}
corresponding to the totally skew-symmetric representations of
$O_{\ell+1}$. The unique   vector
$\phi_{\epsilon}\in\CV^W_{\epsilon,\gamma}$ normalized by $\phi_{\e}(1)=1$ and satisfying the relation
 \be\label{AdM}
  \phi_{\epsilon}(mkm^{-1})\,=\,\phi_{\epsilon}(k)\,,\qquad m\in M\,,
 \ee
is called the $\epsilon$-spherical vector in
$\CV_{\epsilon,\gamma}$. Explicitly, the $\epsilon$-spherical vector
$\phi_{\epsilon}\in\CV^W_{\epsilon,\gamma}$ allows the following
matrix element presentation
 \be\label{genspher1}
  \phi_{\epsilon}(g)\,
  =\,\bigl(v_{\epsilon},\,\pi_W(k)\,v_{\epsilon}\bigr)\,
  \,\chi^B_{\gamma}(a), \qquad g=kan\,,
 \ee
provided by the Iwasawa decomposition \eqref{Iww}.
\end{de}
Given  a generalized spherical vector
$\phi_{\e}\in\CV_{\e,\g}$, we also define the generalized spherical and
Whittaker functions by the following matrix elements:
 \be\label{SPHERF}
 \Phi_{\epsilon,\gamma}(g)=\<\phi_\epsilon,\,
 \pi_{\epsilon,\gamma}(g)\,\phi_\epsilon\>\,,
 \ee
 \be\label{WHITF}
  \Psi_{\epsilon,\gamma}(g)=\<\phi_\epsilon,\,\pi_{\epsilon,\gamma}(g)\,\psi\>\,.
  \ee
The matrix elements above allow explicit integral representations in
terms of the $GL_{\ell+1}(\IR)$-invariant pairing in
$\CV_{\epsilon,\gamma}$. In particular, using the spherical model
\eqref{ModelK} for $\CV_{\epsilon,\gamma}$ the  generalized
spherical function \eqref{SPHERF} may be written as follows:
 \be\label{SPHERF1}
  \Phi_{\epsilon,\gamma}(g)=
  \<\phi_\epsilon,\,\pi_{\epsilon,\gamma}(g)\,\phi_\epsilon\>\,
  =\int\limits_{O_{\ell+1}}\!dk\,\,
  \overline{\phi_\epsilon(k)}\, \phi_{\epsilon}(g^{-1}k)\,.
 \ee
 Then we have
 \be
  \Phi_{\epsilon,\gamma}(1)\,
  =\,d_{|\epsilon|}^{-1}\,,\quad
  d_{|\epsilon|}=\dim \,W_{|\e|}
  =\dim\,\wedge^{|\e|}\IC^{\ell+1}
  =\frac{(\ell+1)!}{|\e|!\,(\ell+1-|\e|)!}\,.
 \ee
as a consequence of  the orthogonality relations for matrix elements of
irreducible $O_{\ell+1}$-representations
$(\pi_i,V_{\mu_i})\in\Irr(O_{\ell+1}),\,i=1,2$:
 \be\label{ORTH}
  \int\limits_{O_{\ell+1}}\!dk\,\,(v_1,\,\pi_{\mu_1}(k)\,w_1)\,
  (v_2,\,\pi_{\mu_2}(k^{-1})\,w_2)
  =\,\frac{\delta_{\mu_1,\mu_2}}{\dim V_{\mu_1}}\,
  (v_1,w_2)\,(v_2,w_1)\,,
 \ee
and the fact that the basis vectors $v_{\e}\in W$ in \eqref{Wbasis}
are orthonormal.

\section{The ramified  Hecke algebra}

 In this Section we introduce an extension
$\CH_r(GL_{\ell+1}(\IR),O_{\ell+1})$ of the commutative spherical
Hecke algebra $\CH(GL_{\ell+1}(\IR),O_{\ell+1})$  retaining its
commutativity property. Thus  defined extended algebra  has a
rich structure of one-dimensional representations and the
$\e$-spherical vectors introduced in Definition \ref{DEF}  appear to
be a common eigenfunction of the elements of
$\CH_r(GL_{\ell+1}(\IR),O_{\ell+1})$.

 Recall that by the Peter-Weyl theory, the convolution algebra
$(L^1(K),*)$ of integrable functions on a compact group $K$ contains
a dense subspace spanned by  the matrix elements
 \be
  F_{\mu,v_i,v_j}(k):=(v_i,\,\pi_\mu(k)\,v_j)\,,\qquad\mu\in\Irr(K)\,,\qquad k\in K\,,
 \ee
where $\{v_i\,,\,i=1,\cdots ,\dim V_{\mu}\}$  provides an
orthonormal basis in each $V_{\mu}$. The normalized matrix elements
supplied with  convolution multiplication may be identified with the
matrix units $E_{ij}^{(\mu)} \in {\rm End}(V_{\mu})$
satisfying the standard relations:
 \be
  E^{(\mu_1)}_{i_1j_1}\cdot
  E^{(\mu_2)}_{i_2j_2}=\,
  \delta_{\mu_1,\mu_2}\,\delta_{j_1,i_2}\,\,E^{(\mu_1)}_{i_1j_2}\,,
 \ee
via
 \be \label{Eme}
  \dim(V_{\mu})\,(v_i,\,\pi_\mu(k)\,v_j)\,\longmapsto\,E_{ij}^{(\mu)}\,.
 \ee

\begin{prop}
Let $K$ be a  compact Lie group and let $L^2(K)$ be the Hilbert
space of square integrable functions on $K$ with the Hermitian
scalar product $(\,,\,)$. Consider a linear subspace of $L^2(K)$
spanned by the diagonal matrix elements in all irreducible
representations $(\pi_{\mu},V_{\mu})$ of $K$
 \be \label{DK}
  D(K)=\bigoplus_{\mu\in{\rm Irr}(K)} D_{\mu}(K)\,,\\
  D_{\mu}(K)\,
  =\,\Span\bigl\{F_{\mu,v_i}(k):=(v_i,\,\pi_\mu(k)\,v_i),\,v_i\in V_{\mu}\bigr\}\,,
 \ee
where $\{v_i\,,\,i=1,\cdots ,\dim V_{\mu}\}$  provides an
orthonormal basis in $V_{\mu}$.
  Then $D(K)$ supplied with convolution multiplication
 \be
  \CH^D(K)=(D(K),*)\,,
 \ee
is an associative commutative algebra.
\end{prop}

 \proof The assertion follows from the identification \eqref{Eme}, but also may be
checked directly   as follows.
 For any pair $(\pi_{\mu_i},V_{\mu_i}),\,i=1,2$ of
$K$-irreducible representations, the convolution of the diagonal
matrix elements $F_{\mu_1,v}$ and $F_{\mu_2,v'}$ appearing in
\eqref{DK}  ($v\in V_{\mu_1}$ and $v'\in V_{\mu_2}$ being  elements of the orthogonal bases)   reads
 \be
  (F_{\mu_1,v}*F_{\mu_2,v'})(\tilde{k})\,
  =\int\limits_K\!dk\,F_{\mu_1,v}(k^{-1})\,F_{\mu_2,v'}(k\tilde{k})\\
  =\int\limits_K\!dk\,\ov{(v,\,\pi_{\mu_1}(k)\,v)}\,\,
  (v',\,\pi_{\mu_2}(k)\,\pi_{\mu_2}(\tilde{k})\,v')\\
  =\,\delta_{\mu_1,\mu_2}\,\frac{(v,\,v')}{\dim V_{\mu_1}}\,
  (v,\,\pi_{\mu_1}(\tilde{k})\,v')\,.
 \ee
The latter equality follows by the orthogonality relations
\eqref{ORTH}. Taking into account the condition
 $(v,\,v')=\delta_{v,v'}$ we obtain
the commutativity of $(\CH^D(K),*)$. $\Box$

There are various finite-dimensional sub-algebras of $\CH^D(K)$
obtained by restricting the direct sum  \eqref{DK} to a finite
subset of ${\rm Irr}(K)$. In the following we consider the case of
$K=O_{\ell+1}$ and the subset  ${\rm Irr}_W(O_{\ell+1})$ of the
irreducible totally skew-symmetric
$O_{\ell+1}$-representations occurring in the decomposition of
$(\pi_W,W)$,
 \be
  W\,=\,\bigoplus_{k=0}^{\ell+1}\,
  W_k=\bigoplus_{k=0}^{\ell+1} \wedge^k \IC^{\ell+1}, \qquad
  d_k=\dim \,W_k=\frac{(\ell+1)!}{k!(\ell+1-k)!}\,.
  \ee
Denote by $F_{v_{\e}}$ the corresponding diagonal  matrix elements
in the basis $\{v_\epsilon:\,\e\in\CP\}\subset W$ introduced in
\eqref{Wbasis} (the specific $O_{\ell+1}$-irreducible representation
is fixed by the value $|\e|$). Note that the functions $F_{v_{\e}}$ satisfy
the condition \eqref{AdM} of invariance under the adjoint action of
the subgroup $M\subset O_{\ell+1}$.

\begin{lem} The following element of $\CH^D(O_{\ell+1})$,
 \be\label{Proj}
  \Delta_W(k)\,
  =\sum_{\epsilon\in\CP} d_{|\epsilon|}
  \<v_{\epsilon},\pi_W(k)\,v_{\epsilon}\>\,,\qquad k\in
  O_{\ell+1}\,,
  \ee
  acting in $L^2(O_{\ell+1})$ via convolution is a projector
 \be
 \Delta_W*\Delta_W=\Delta_W\,,
 \ee
 onto a subspace isotypical to $W$ i.e.   allowing a decomposition into the
same set of irreducible representations of $O_{\ell+1}$
but possibly with different multiplicities.
Under isomorphism \eqref{Eme} the  projector
 is given by
\be
 \Delta_W=\sum_{k=0}^{\ell+1}  \,\,{\rm Id}_{W_k}=
\sum_{k=0}^{\ell+1}  \,\sum_{j=1}^{\dim W_k} \,E^{(\pi_k)}_{ii}\,.
\ee
\end{lem}
\proof Direct verification using the orthogonality relations
\eqref{ORTH}. $\Box$

Now using the properties of $\CH^D(O_{\ell+1})$ we construct an
extension of the spherical Hecke algebra
$\CH(GL_{\ell+1}(\IR),O_{\ell+1})$ of $O_{\ell+1}$-bi-invariant
functions on $GL_{\ell+1}(\IR)$. We start with a construction of the
following non-commutative  extension of
$\CH(GL_{\ell+1},O_{\ell+1})$. Consider the  space ${\rm
Fun}(GL_{\ell+1}(\IR))$ as a module under right and left actions of
the orthogonal subgroup $O_{\ell+1}$. Introduce
a subspace $\Fun_W(GL_{\ell+1}(\IR))\subset\Fun_(GL_{\ell+1}(\IR))$  of functions
transforming  under the right and  left  $O_{\ell+1}$-actions via
the representations  $W$ and its dual $W^*$, correspondingly.
The linear space ${\rm Fun}_W(GL_{\ell+1}(\IR))$ has a natural structure
of a module over the algebra of  $O_{\ell+1}$-bi-invariant functions
on $GL_{\ell+1}$. Explicitly,
the following set of functions provides a basis in  ${\rm
Fun}_W(GL_{\ell+1}(\IR))$
 \be\label{gen1}
  F_{v_{\e},v_{\e'}}(g)\,
  =\,(v_\e, \pi_W(g),v_{\e'})\,F(g)\,,\qquad F(k_1gk_2)=F(g), \quad
  k_1,k_2\in O_{\ell+1}\\
  F_{v_{\e},v_{\e'}}(g)=0, \qquad |\epsilon|\neq |\epsilon'|,
 \ee
where the representation $(\pi_W,W)$ of $GL_{\ell+1}(\IR)$,
the pairing $(\,,\,)$ and  the orthogonal bases
 $\{v_\e,\,\e\in \CP\}$ are defined in \eqref{Fundrep},
 \eqref{Wbasis}. Note that the subspace ${\rm Fun}_W(GL_{\ell+1}(\IR))$ may be defined as
the subspace of ${\rm Fun}(GL_{\ell+1}(\IR))$ such that under the
left and right action the projector $\Pi_W$ acts on this subspace by the identity operator.

\begin{lem} \label{Conm}
  Consider the functions on $GL_{\ell+1}(\IR)$  of the form \eqref{gen1}:
 \be\label{FF}
  F_{v_\e, v_{\e'}}(g)=(v_\epsilon,\,\pi_W(g)\,v_{\e'})\,F(g)\,,\quad\e,\e'\in\CP\,,
 \ee
where $F\in \CH(GL_{\ell+1}(\IR),O_{\ell+1})$ i.e. $O_{\ell+1}$-bi-invariant function on $GL_{\ell+1}(\IR)$,
 \be
  F(k_1gk_2)\,=\,F(g),\qquad k_1,k_2\in O_{\ell+1}\,,
 \ee
satisfying the conditions of rapid decay at infinity, and
$\{v_{\epsilon},\,\e\in\CP\}$ is the basis \eqref{Wbasis} of the
representation $W=\oplus_k \wedge^k \IC^{\ell+1}$. Then for the
convolution of two such functions,
 \be
  F_{v_{\e_1}, v_{\e_1'}}(g)=(v_{\epsilon_1},\,\pi_W(g)\,v_{\epsilon'_1})\,F(g)\,,\\
  G_{v_{\e_2},v_{\e_2'}}(g)=(v_{\epsilon_2},\,\pi_W(g)\,v_{\epsilon'_2})\,G(g)\,,
 \ee
the following relation holds
 \be
  (F_{v_{\e_1}, v_{\e_1'}}*G_{v_{\e_2}, v_{\e_2'}})(g)\,
  =\frac{\delta_{\epsilon'_1,\epsilon_2}}{\dim W_{|\epsilon_1'|}}\,
  (v_{\epsilon_1},\,\pi_W(g)\,v_{\epsilon_2'})\,(F*G)(g)\,.
 \ee
\end{lem}

\proof In straightforward way, we have
\be
 (F_{v_{\e_1}, v_{\e_1'}}*G_{v_{\e_2}, v_{\e_2'}})(\tilde{g})\,
  =\!\!\int\limits_{GL_{\ell+1}(\IR)}\!\!\!d\mu(g)\,\,
  F_{v_{\e_1},v_{\e_1'}}(g^{-1})\,\,G_{v_{\e_2}, v_{\e_2'}}(g\tilde{g})\\
  =\!\!\int\limits_{GL_{\ell+1}(\IR)}\!\!\!d\mu(g)\,\,
  (v_{\epsilon_1},\,\pi_W(g^{-1})\,v_{\epsilon'_1})\,\,
  (v_{\epsilon_2},\,\pi_W(g\tilde{g})\,v_{\epsilon'_2})\,\,
  F(g^{-1})\,G(g\tilde{g})\,.
 \ee
Applying the Iwasawa  decomposition \eqref{Iww} we obtain
 \be
  (F_{v_{\e_1}, v_{\e_1'}}*G_{v_{\e_2}, v_{\e_2'}})(\tilde{g})
  =\int\limits_{O_{\ell+1}\times A \times N}\!\!dk\,da\,dn\,e^{-2\rho(a)}\,\,
  F(n^{-1}a^{-1})\,G(an\tilde{g})\\
 \times \,\,
  (v_{\epsilon_1},\,\pi_W(n^{-1}a^{-1}k^{-1})\,v_{\epsilon'_1})\,
  (v_{\epsilon_2},\,\pi_W(kan \tilde{g})\,v_{\epsilon'_2})\,.
 \ee
We calculate the integral over $k$ using the orthogonality
relations \eqref{ORTH}:
 \be
  \int\limits_{O_{\ell+1}}\!\!dk\,\,
  (v_{\epsilon_1},\,\pi_W(n^{-1}a^{-1}k^{-1})\,v_{\epsilon'_1})\,
  (v_{\epsilon_2},\,\pi_W(kan \tilde{g})\,v_{\epsilon'_2})\\
  =\,\frac{\delta_{\e_1',\e_2}}{d_{|\e_1'|}}\,
  (v_{\epsilon_1},\,\pi_W(\tilde{g})\,v_{\epsilon_2'})\,.
 \ee
Hence we obtain the following:
 \be
  (F_{v_{\e_1}, v_{\e_1'}}*G_{v_{\e_2}, v_{\e_2'}})(\tilde{g})\\
  =\,\frac{\delta_{\e'_1,\e_2}}{d_{|\e'_1|}}\,
  (v_{\epsilon_1},\,\pi_W(\tilde{g})\,v_{\epsilon'_2})\!
  \int\limits_{A \times N}\!\!da\,dn\,e^{-2\rho(a)}\,\,
  F(n^{-1}a^{-1})\,G(an\tilde{g})\\
  =\,\frac{\delta_{\e'_1,\e_2}}{d_{|\e'_1|}}\,
  (v_{\epsilon_1},\,\pi_W(\tilde{g})\,v_{\epsilon'_2})\!
  \int\limits_{GL_{\ell+1}}\!\!d\mu(g)\,\,
  F(g^{-1})\,G(g\tilde{g})\\
  =\,\frac{\delta_{\e'_1,\e_2}}{d_{|\e'_1|}}\,
  (v_{\epsilon_1},\,\pi_W(\tilde{g})\,v_{\epsilon'_2})\,(F*G)(\tilde{g})\,,
 \ee
which completes our proof. $\Box$

From  Lemma  \ref{Conm} we infer that the linear subspace ${\rm Fun}_W(GL_{\ell+1}(\IR))
 \subset {\rm Fun}(GL_{\ell+1}(\IR)) $
is closed under convolution and  is  isomorphic to the algebra
of $O_{\ell+1}$-bi-invaraint functions
 \be
  \CH_W(GL_{\ell+1}(\IR),O_{\ell+1})
  =({\rm Fun}(GL_{\ell+1}(\IR))\otimes\End_{gr}(W))^{O_{\ell+1}\times O_{\ell+1}}\,,
 \ee
 where
 \be
  \End_{gr}(W)\,=\,\bigoplus_{k=0}^{\ell+1}\,\,{\rm End}(W_k)\,,
  \ee
so that under  left and right actions  $O_{\ell+1}$ acts diagonally on
both factors.

We are interested in the maximal commutative subalgebra of
$\CH_W(GL_{\ell+1}(\IR))$, called ramified Hecke algebra
$\CH_r(GL_{\ell+1},O_{\ell+1})$.

\begin{prop} \label{COR1}
The $M$-invariant linear subspace of
$\Fun_W(GL_{\ell+1}(\IR))$ consisting of the functions \eqref{FF} satisfying
 \be\label{Minv2}
  F_{v_{\e},v_{\e'}}(mgm^{-1})=F_{v_{\e},v_{ \e'}}(g), \qquad m\in M\,,
 \ee
is spanned by the following functions, for $\e\in\CP$,
 \be\label{F}
  F_{v_{\e}}(g)=(v_\epsilon,\,\pi_W(g)\,v_{\e})\,F(g)\,,\quad F(k_1gk_2)=F(g)\,,
 \ee
and provides  a  maximal commutative subalgebra of $\CH_W(GL_{\ell+1},O_{\ell+1})$.
Moreover, the  following relations for \eqref{F} hold
 \be
  F_{v_{\e_1}}*G_{v_{e_2}}\,
  =\,\frac{\delta_{\e_1, \e_2}}{\dim W_{|\e_1|}}\,
  (F*G)_{v_{\e_1}}\,.
 \ee
\end{prop}
\proof Directly follows by Lemma \ref{Conm}. $\Box$

\begin{lem} As a linear space, $\CH_r(GL_{\ell+1},O_{\ell+1})$ is
isomorphic to the space  of $\CW$-invariant
  functions on the centralizer $A\times M$ of $A$,
  \be
 \CH_r(GL_{\ell+1},O_{\ell+1})\simeq {\rm Fun}(A\cdot M)^{\CW}\,.
\ee
where $\CW\simeq\mathfrak{S}_{\ell+1}$ is the Weyl group of $(GL_{\ell+1}(\IR),\,A)$.
\end{lem}

\proof Indeed we have
 \be
  \CH_r(GL_{\ell+1}(\IR),O_{\ell+1})\,
  =\,\Big(\Fun(GL_{\ell+1}(\IR))\otimes(\End_{gr}(W))^M\Big)^{O_{\ell+1}\times O_{\ell+1}}\,.
 \ee
Consider the polar Cartan covering map
 \be
 O_{\ell+1}\times A\times  O_{\ell+1}\longrightarrow GL_{\ell+1}(\IR)=
O_{\ell+1}\, A\,  O_{\ell+1}\,,
\ee
with a fiber isomorphic to $M\rtimes \CW$, then we deduce
\be
 \CH_r(GL_{\ell+1},O_{\ell+1})\,
 =\,\Big(\Fun(A) \otimes \End_{gr}(W)^M\Big)^\CW\,.
\ee
Taking into account the isomorphisms
 \be
  (\End_{gr}(W)^{M}\simeq \IC[M]\,,\qquad
\Fun(A\times M)=\Fun(A)\otimes \Fun(M)\,,
\ee
we arrive at the required statement. $\Box$

Thus Proposition  \ref{COR1} provides a construction of an extension
$\CH_r(GL_{\ell+1}(\IR),O_{\ell+1})$ of the spherical Hecke algebra
$\CH(GL_{\ell+1}(\IR),O_{\ell+1})$.

\begin{de} \label{RHA} The commutative ramified Hecke algebra
$\CH_r(GL_{\ell+1}(\IR),O_{\ell+1})$ is defined as the convolution
algebra of the   elements of the linear subspace of ${\rm Fun}_W(GL_{\ell+1}(\IR))$
invariant under the adjoint action of $M$ via \eqref{Minv2}.
This subspace is spanned by the  functions \eqref{F}:
 \be
  F_{\epsilon}(g)\,
  =\,(v_\epsilon, \pi_W(g)\,v_\epsilon)\,F(g), \qquad
  F(g)\in \CH(GL_{\ell+1}(\IR),O_{\ell+1})\,,\qquad \epsilon\in \CP\,.
  \ee
  \end{de}

In the next Section we define a one-parameter family of elements in the
ramified Hecke algebra $\CH_r(GL_{\ell+1}(\IR),O_{\ell+1})$
generalizing the spherical Hecke-Baxter operator \eqref{HB0}.

\section{The $GL_{\ell+1}(\IR)$ Hecke-Baxter operator in ramified case}

In this Section we define the  Hecke-Baxter operator acting on generalized
spherical vectors $\phi_{\e}\in\CV_{\e,\g}$ proposed in Definition \ref{DEF}
as a one-parameter family of elements of the
ramified Hecke algebra $\CH_r(GL_{\ell+1}(\IR),O_{\ell+1})$
introduced in Definition \ref{RHA}. The main requirement is that
  its action on the generalized spherical vectors in the
principal series representation $\CV_{\epsilon,\gamma}$ (and thus on
  the corresponding generalized spherical and
Whittaker functions \eqref{SPHERF},\eqref{WHITF}) should be given by
multiplication by the corresponding  $L$-factor \eqref{Lf}.

\begin{te}\label{MAINTHM}
The  Hecke-Baxter operator defined as convolution with the
one-parameter family of functions
 \be\label{GenQ123}
  \wh{Q}_s(g)\,
  =\,\Delta_W(g)\,Q_s(g)
  =\,\Delta_W(g)\,(c\pi^{-1})^{\frac{\ell(\ell+1)}{4}}\,\,|\det g|^{s+\frac{\ell}{2}}\,
  e^{-c\Tr(g^{\top}g)},
 \ee
with
 \be
  \Delta_W(g)\,
  =\sum_{\epsilon\in\CP} d_{|\epsilon|}\,
  (v_\epsilon,\pi_W(g)\,v_{\epsilon})\,,\qquad d_{|\e|}=\dim W_{|\e|},
  \qquad W_{|\e|}=\wedge^{|\e|}\,\IC^{\ell+1}\,,
 \ee
is acting on the $\epsilon$-spherical vector \eqref{genspher1},
 \be
  \phi_{\epsilon}(g)\,
  =\,(v_{\epsilon},\pi_W(k(g))\,v_{\epsilon})\,\,
  \chi^B_{\gamma}(a),\qquad
  g=kan\in O_{\ell+1}AN\,,
 \ee
via multiplication by the Archimedean $L$-factor \eqref{Lf},
 \be\label{eL1}
  L(s|\epsilon,\gamma)\,
  =\prod_{j=1}^{\ell+1}c^{-\frac{s+\epsilon_j-\imath \gamma_j}{2}}\,
  \Gamma\left(\frac{s+\epsilon_j-\imath \gamma_j}{2}\right)\,.
 \ee
\end{te}
\proof The proof is given in Appendix.

\begin{cor}
The Hecke-Baxter operator defined as convolution with one-parameter
family of functions \eqref{GenQ123}
 \be\label{GenQ123a}
  \wh{Q}_s(g)\,
  =\,\Delta_W(g)\,(c\pi^{-1})^{\frac{\ell(\ell+1)}{4}}\,\,|\det g|^{s+\frac{\ell}{2}}\,
  e^{-c\Tr(g^{\top}g)},
 \ee
acts on the generalized spherical and Whittaker functions
\eqref{SPHERF} and \eqref{WHITF},
 \be
  \Phi_{\epsilon,\gamma}(g)=\<\phi_\epsilon,\pi_{\epsilon,\gamma}(g)
  \,\phi_\epsilon\>, \qquad
  \Psi_{\epsilon,\gamma}(g)=\<\phi_\epsilon,\pi_{\epsilon,\gamma}(g)
  \,\psi\>\,,
  \ee
via multiplication by the  Archimedean $L$-factor
 \be
  L(s|\epsilon,\gamma)\,
  =\prod_{j=1}^{\ell+1}c^{-\frac{s+\epsilon_j-\imath \gamma_j}{2}}\,
  \Gamma\left(\frac{s+\epsilon_j-\imath \gamma_j}{2}\right)\,.
 \ee
\end{cor}

Thus the  generalized  Hecke-Baxter operator  is obtained from the spherical
 Hecke-Baxter operator \eqref{HB0},
 \be\label{GenQ1234}
  Q_s(g)\,
  =\,(c\pi^{-1})^{\frac{\ell(\ell+1)}{4}}\,|\det g|^{s+\frac{\ell}{2}}\,
  e^{-c\Tr(g^{\top}g)},
 \ee
via adding the prefactor  $\Delta_W(g)$
 \be
  \Delta_W(g)\,
  =\sum_{\epsilon\in\CP}
  d_{|\e|}\,(v_{\epsilon},\pi_W(g)\,v_{\epsilon})\,.
  \ee
This prefactor allows various useful interpretations.
First thing to notice is that $\Delta_W(g)$ is a lifting to
$GL_{\ell+1}(\IR)$ of the projector \eqref{Proj} onto the
$O_{\ell+1}$-subrepresentation $W\subset L^2(O_{\ell+1})$. Let us
 use the natural grading by $k$ on the  vector space $W=\oplus_k
\wedge^k \IC^{\ell+1}$ to introduce the  vector space of
endomorphisms of $W$ respecting this grading
 \be
  \End_{gr}(W)=\bigoplus_{k=0}^{\ell+1}\End(W_k)\,.
 \ee
This space has a natural structure of a
$GL_{\ell+1}(\IR)\times GL_{\ell+1}(\IR)$-module under the left and
right actions of $GL_{\ell+1}(\IR)$, hence we have following representation
of $\Delta_W(g)$:
\be
\Delta_W(g)\,=\,\Tr_{\End_{gr}(W)}(1\otimes g)\,.
\ee

Now, let us consider the universal enveloping algebra
$\CU\mathfrak{gl}_{\ell+1}$ of the Lie algebra
$\mathfrak{gl}_{\ell+1}=\Lie(GL_{\ell+1}(\IR))$. Let
$\pi_W(\CU\mathfrak{gl}_{\ell+1})\subset\End(W)$ be the image of
$\CU\mathfrak{gl}_{\ell+1}$ given by the
$\CU\mathfrak{gl}_{\ell+1}$-representation $(\pi_W,W)$.
The algebra $\pi_W(\CU\mathfrak{gl}_{\ell+1})$ may be succinctly
described  in terms of the Clifford algebra ${\rm
Cliff}_{2\ell+2}={\rm Cliff}(\IC^{2\ell+2})$ generated by
$\psi^*_i$, $\psi_i$, $i=1,\ldots ,(\ell+1)$ subjected to the following
anti-commutator relations:
 \be
  [\psi_i,\psi_j]_+=0, \quad [\psi^*_i,\psi^*_j]_+=0, \quad
  [\psi^*_i,\psi_j]_+=\delta_{ij}, \quad i,j=1,\ldots ,(\ell+1)\,.
 \ee
 The Clifford algebra has a structure of a graded algebra upon assigning the
following grading to the generators:   ${\rm deg}(\psi_i^*)=1$,  ${\rm
  deg}(\psi_i)=-1$. This grading is compatible with the   grading in
$W$ considered as a representation of ${\rm Cliff}_{2\ell+2}$. Let
${\rm Cliff}^{(0)}_{2\ell+2}\subset {\rm Cliff}_{2\ell+2}$ be the
zero grade sub-algebra. Then we  have  the following identification:
 \be
  \pi_W(\CU\mathfrak{gl}_{\ell+1})\simeq   {\rm Cliff}^{(0)}_{2\ell+2}
  \simeq {\rm End}_{gr}(W) \,,
 \ee
and therefore we arrive at the following presentation of $\Delta_W(g)$
 \be
  \Delta_W(g)\,
  =\,\Tr_{{\rm Cliff}_{2\ell+2}^{(0)}}(1\otimes g)\,,
 \ee
which might be considered as an analog of the classical Cauchy
identity for the associative algebra ${\rm Cliff}_{2\ell+2}^{(0)}$.

\section{Fixing ambiguities via Fourier  transform }

In the main part of this note we have considered various
representation theory constructions associated to Lie group
$GL_{\ell+1}(\IR)$. There is a kind  of intrinsic arbitrariness in
these constructions manifested  in particular in the existence of
the parameter $c\in\IR^*_+$ entering expressions of the local Archimedean
$L$-factors \eqref{Lf},\eqref{Lf0} and of the Hecke-Baxter operators
\eqref{HB0},\eqref{GenQ123}. One might fix this ambiguity by  taking into
account the considerations related with the theory of global
zeta-functions encompassing the information about different
completions of rational numbers (ultimately this reduces to invoking
integral structure on the reals). Precisely we  require the
existence of  simple functional equations for the global zeta functions
constructed as a product of local factors. In terms of the local
Archimedean $L$-factors this, in particular, favours  simple transformation
properties for the Gaussian measure,
 \be\label{Gmes}
  d\g(g)=e^{-c\Tr(g^{\top}g)}\,dg\,,\qquad
  dg\,=\prod_{i,j=1}^{\ell+1}dg_{ij}\,,
 \ee
entering the integral expression for the local Archiemdean
$L$-factor. For instance, local Archimedean $L$-factor \eqref{Lf0}
for the spherical principal series representation
$(\pi_{0,0},\CV_{0,0})$ induced via  the $B$-character
$\chi^B_{0}$ reads from \eqref{EIGENVAL}:
 \be
  L(s|0)\,=\,(c\pi^{-1})^{\frac{\ell(\ell+1)}{4}}
  \!\int\limits_{GL_{\ell+1}(\IR)}\!\!d\mu(g)
  \,\,|\det g|^{s+\frac{\ell}{2}}\, e^{-c\Tr(g^{\top}g)}\\
  =\,(c\pi^{-1})^{\frac{\ell(\ell+1)}{4}}\!\int\limits_{GL_{\ell+1}(\IR)}\!\!d\g(g)\,\,
  |\det g|^{s-\frac{\ell}{2}-1}\,,
 \ee
where $d\mu(g)$ is the Haar measure \eqref{Hmes} on $GL_{\ell+1}(\IR)$:
 \be
  d\mu(g)\,=\,|\det(g)|^{-(\ell+1)}\,dg\,.
 \ee
Consider the Fourier transform on $\Mat_{\ell+1}(\IR)$:
 \be\label{FT1}
  (\CF f)(\tilde{g})\,
  =\!\!\int\limits_{\Mat_{\ell+1}(\IR)}\!\!
  dg\,\,\,e^{2\pi\imath\,\Tr(g^{\top}\tilde{g})}\,f(g)\,.
 \ee
This is a unitary operator on
$L^2\bigl(\Mat_{\ell+1}(\IR),dg\bigr)$.

 \begin{lem} The Gaussian measure \eqref{Gmes} with $c=\pi$,
 \be
  d\g(g)\,
  =\,G(g)\,dg\,
  =\,e^{-\pi\Tr(g^{\top}g)}\,dg\,,
 \ee
is self-dual with respect to the Fourier transform \eqref{FT1}:
 \be
  (\CF G)(g)\,=\,G(g)\,.
 \ee
\end{lem}
\proof The assertion directly follows from the standard integral
identity:
 \be\label{FGauss}
  \int\limits_{\IR}\!dx\,\,
  e^{2\pi\imath yx\,-\,\pi x^2}\,
  =\,e^{-\pi y^2}\,.
 \ee
$\Box$

Therefore, it seems  natural to specify $c=\pi$ in expressions for
the local Archimedean $L$-factors \eqref{Lf0}.  However,
the proposed construction of the non-spherical Hecke-Baxter operator
implies another choice.  By Theorem \ref{MAINTHM}, the non-spherical Hecke-Baxter operator is given by the convolution with the function \eqref{GenQ123},
 \be\label{GenQ11}
  \wh{Q}_s(g)\,
  =\,\Delta_W(g)\,(c\pi^{-1})^{\frac{\ell(\ell+1)}{4}}\,|\det g|^{s+\frac{\ell}{2}}\,
  e^{-c\,\Tr(g^{\top}g)}\,,
 \ee
where
 \be
  \Delta_W(g)\,
  =\sum_{\epsilon\in\CP} d_{|\epsilon|}\,
  (v_\epsilon,\,\pi_W(g)\,v_{\epsilon})\,,
 \ee
hence a natural requirement to impose would be to have  simple
transformations properties with respect to
Fourier transform of the  modified Gaussian measure
 \be\label{ModGauss}
  d\g_W(g)\,
  =\,G_W(g)\,dg\,
  =\,\Delta_W(g)\,(c\pi^{-1})^{\frac{\ell(\ell+1)}{4}}\,e^{-c\,\Tr(g^{\top}g)}\,dg\,.
 \ee
The choice of $c=\pi$ leads however to a non-trivial
transformation of the measure \eqref{ModGauss}.

\begin{lem} The modified  Gaussian measure \eqref{ModGauss} with $c=\pi$
 \be
  d\g_W(g)\,=\,G_W(g)\,dg\,
  =\,\Delta_W(g)\,e^{-\pi\Tr(g^{\top}g)}\,dg\,,
 \ee
satisfies the following identity
 \be\label{IDFT}
  (\CF G_W)(g)\,=\,\tilde{G}_W(g)\,,\\
  \tilde{G}_W(g)\,
  =\,\Delta_W(\imath g)\,(c\pi^{-1})^{\frac{\ell(\ell+1)}{4}}\,e^{-\pi\Tr(g^{\top}g)}\,.
 \ee
\end{lem}
\proof The assertion follows by the simple  identities:
 \be\label{FT124}
  \CF(G)(x)=G(x)\,, \qquad G(x)=e^{-\pi x^2}\,,\\
  \CF(G^{(1)})(x)=\imath G^{(1)}(x), \qquad G^{(1)}(x)=xe^{-\pi x^2}\,,
 \ee
where the later identity easily follows from the former one  via
 \be
  G^{(1)}(x)=\frac{1}{2\pi\imath}\frac{\pr}{\pr x} G(x)\,.
 \ee
Precisely, we should take  into  account that $\Delta_W(g)$
 \be\label{DeltaW}
  \Delta_W(g)=\sum_{I_k} \,d_k\,\,\Delta_{I_k}(g)\,,\quad
  I_k=(i_1<\cdots <i_k),\,\,1\leq i_a\leq (\ell+1)\,,
 \ee
is a sum of principal $k$-minors $\Delta_{I_k}(g)$. Hence
the r.h.s. of \eqref{DeltaW} is a sum of monomials  such that
each element $g_{ij}$ enters the product in a power at most one.
Therefore, applying the integral identity\eqref{FT124} to each
variable $g_{ij}$ yields the identity \eqref{IDFT}. $\Box$

Thus, the modified Gauss measure $d\g_W(g)$ transforms non-trivially under Fourier
transform. To resolve this conundrum  one should recall that  the requirement for
the real measure $d\g_W(g)$ to be self-dual with respect to Fourier transform may be
weaken to the condition on a complex measure to be self-dual under a combination of the
Fourier transform and complex conjugation. This leads to considering
the following Feynmann type measure (understood as usual as a
limit of  a well-defined complex measure).

\begin{prop}
The modified  Feynmann (imaginary Gaussian)   measure
 \be
  d\,\wt{\g}_W(g)\,
  =\,\wt{G}_W(g)\,dg\,
  =\,\Delta_W(g)\,(c\pi^{-1})^{\frac{\ell(\ell+1)}{4}}\,e^{-\imath\pi\,\Tr(g^{\top}g)}\,dg\,,
 \ee
satisfies the following identity
 \be\label{IDFT1}
  (\CF\,\wt{G}_W)(g)\,
  =\,e^{-\frac{\imath\pi (\ell+1)^2}{4}}\,\,\overline{\wt{G}_W(g)}\,.
 \ee
\end{prop}
\proof This  directly follows from the previous Lemma and
the following simple computations:
 \be\label{FT1245}
  (\CF\,\wt{G})(x)=e^{-\frac{\imath\pi}{4}}\,\overline{\wt{G}(x)}, \qquad
  \wt{G}(x)=e^{-\imath \pi x^2},\\
  (\CF\,\wt{G}^{(1)})(x)
  =e^{-\frac{\imath\pi}{4}}\, \overline{\wt{G}^{(1)}(x)}, \qquad
  \wt{G}^{(1)}(x)=xe^{-\imath \pi x^2}\,.
 \ee
To calculate the Fourier transform  we use deformations of the
basic functions by $\varepsilon>0$:
 \be
  \wt{G}(g|\varepsilon)=e^{-(\varepsilon+\imath \pi)x^2}, \qquad
  \wt{G}^{(1)}(g|\varepsilon)=xe^{-(\varepsilon+\imath \pi)x^2}, \qquad
  \varepsilon\to 0_+\,,
 \ee
to render the integrals well-defined. $\Box$

Arithmetic implications of replacing the
Gaussian  quadratic measure by the Feynmann one is an interesting direction to pursue.

\section{Appendix: Proof of Theorem \ref{MAINTHM}}

Introduce the following notation for the components of an element
$g\in GL_{\ell+1}(\IR)$ upon the Iwasawa decomposition \eqref{Iwasawa}:
\be
g\,=\,k(g)\,a(g)\,n(g)\,, \quad k(g)\in O_{\ell+1}, \quad a(g)\in A, \quad
n(g)\in N\,.
\ee
The convolution action of the function \eqref{GenQ1234},
 \be\label{Q}
  \wh{Q}_s(g)\,
  =\,\sum_{\epsilon'\in\CP} d_{|\epsilon'|}\,
  (v_{\e'},\,\pi_W(g)\,v_{\e'})\,\,
  Q_s(g),\quad
  d_{|\e'|}=\dim\,W_{|\e'|}\,,
 \ee
on the $\e$-spherical vector \eqref{genspher1},
 \be\label{phi}
  \phi_{\epsilon}(g)\,
  =\,(v_{\epsilon},\,\pi_W(k(g))\,v_{\epsilon})\,\,
  \chi^B_{\gamma}(a(g))\,,
 \ee
reads
 \be\label{Qphi0}
  \tilde{\phi}_\epsilon(\tilde{g}):
  =(\wh{Q}_s*\phi_\epsilon)(\tilde{g})\,
  =\!\int\limits_{GL_{\ell+1}(\IR)}\!\!\!d\mu(g)\,\,\,
  \wh{Q}_s(g)\,\,\phi_\epsilon(g^{-1}\tilde{g})\,.
 \ee
Consider the polar Cartan covering map
 \be\label{CPD3}
  O_{\ell+1}\times A\times O_{\ell+1}\longrightarrow
  GL_{\ell+1}(\IR)\,,
 \ee
given by
 \be
  (k_1,a,k_2)\,\longmapsto\,k_1ak_2, \qquad a\in A,\qquad k_1,k_2\in O_{\ell+1}\,.
 \ee
A fiber of this map may be identified with finite group $\CW^M:=M\rtimes \CW$,
where $\CW=\frak{S}_{\ell+1}$ is the Weyl group of $(GL_{\ell+1}(\IR),A)$.
The lift of the Haar measure on $GL_{\ell+1}(\IR)$ to the product of
Lie groups in \eqref{CPD3}  is given by
 \be
  d\mu(g)\,=\,|\Delta(a)|\,dk_1\,da\,dk_2\,,\quad
  \Delta(a)\,=\,\prod_{i<j}\Big(\frac{a_i}{a_j}-\frac{a_j}{a_i}\Big),
  \qquad g=k_1ak_2\,.
 \ee
Therefore, substituting \eqref{Q}, \eqref{phi} into \eqref{Qphi0} we obtain
 \be\label{Qphi1}
  \tilde{\phi}_\epsilon(\tilde{g})\,
  =\frac{1}{|\CW^M|}\,\,\int\limits_{O_{\ell+1}\times A \times
    O_{\ell+1}}
  \!dk_1\,da\,dk_2\,\,|\Delta(a)|\,\,
  Q_s(a)\\
  \times\sum_{\epsilon'\in\CP}d_{|\epsilon'|}\,
  (v_{\epsilon'},\,\pi_W(k_1ak_2)\,v_{\epsilon'})\,
  (v_\epsilon,\,\pi_W\bigl(k(k_2^{-1}(k_1a)^{-1}\tilde{g})\bigr)\,v_\epsilon)\,
  \chi^B_{\gamma}\bigl(a((k_1a)^{-1}\tilde{g})\bigr)\,,
 \ee
 where we use
 \be
 a(k_2^{-1}(k_1a)^{-1}\tilde{g})=a((k_1a)^{-1}\tilde{g})\,,
 \ee
providing that the $\chi^B_{\g}$-factor is left $O_{\ell+1}$-invariant.
Noting that
 \be
  k(k_2^{-1}(k_1a)^{-1}\tilde{g})=k_2^{-1}\,k((k_1a)^{-1}\tilde{g})\,,
 \ee
we calculate the integral over $k_2\in O_{\ell+1}$ applying the
orthogonality relations \eqref{ORTH} as follows:
 \be
  \int\limits_{O_{\ell+1}}\!dk\,
  (v_{\epsilon'},\,\pi_W(k_1ak_2)\,v_{\epsilon'})\,
  (v_\epsilon,\,\pi_W(k_2^{-1})\,
  \pi_W\bigl(k({k_1a}^{-1}\tilde{g}\bigr)\,v_\epsilon)\\
  =\,\frac{\delta_{\epsilon',\epsilon}}{d_{|\epsilon|}}\,
  (v_{\epsilon},\,\pi_W(k_1a)\,\pi_W\bigl(k((k_1a)^{-1}\tilde{g})\bigr)\,v_\epsilon)
  \,.
 \ee
Then \eqref{Qphi1} takes the following form:
 \be\label{Qphi21}
  \tilde{\phi}_\epsilon(\tilde{g})\,
  =\frac{1}{|\CW^M|}\,\,\!\int\limits_{O_{\ell+1}\times A}\!dk_1 da\,|\Delta(a)|\,\,
  Q_s(a)\,\chi^B_{\g}\bigl(a((k_1a)^{-1}\tilde{g})\bigr)\\
  \times\,(v_{\epsilon},\,\pi_W(k_1a)\,
  \pi_W\bigl(k((k_1a)^{-1}\tilde{g})\bigr)\,v_\epsilon)\,.
 \ee
Next, considering the Iwasawa decomposition of the argument
$\tilde{g}=\tilde{k}\tilde{a}\tilde{n}$ we get
 \be
  k((k_1a)^{-1}\tilde{g})=k((k_1a)^{-1}\tilde{k}), \quad
  a((k_1a)^{-1}\tilde{g})= a((k_1a)^{-1}\tilde{k})\tilde{a}\,,
 \ee
which entails
 \be
  \chi^B_{\gamma}\bigl(a((k_1a)^{-1}\tilde{g})\bigr)\,
  =\,\chi^B_{\gamma}(a((k_1a)^{-1}\tilde{k}))\,\,
  \chi^B_{\gamma}(a(\tilde{g}))\,,\\
  (v_{\epsilon},\,\pi_W\bigl(k((k_1a)^{-1}\tilde{g})\bigr)\,v_\epsilon)\,
  =\,(v_{\epsilon},\,\pi_W\bigl(k((k_1a)^{-1}\tilde{k}\bigr)\,v_\epsilon)\,.
 \ee
Hence substitution into \eqref{Qphi1} leads to the following, for
$\tilde{k}=k(\tilde{g})\in O_{\ell+1}$\,,
 \be\label{Qphi2}
  \tilde{\phi}_\epsilon(\tilde{g})\,
  =\,\frac{1}{|\CW^M|}\,\,\chi^B_{\gamma}\bigl(a(\tilde{g})\bigr)\!\!
  \int\limits_{O_{\ell+1}\times A}\!dk_1\,da\,|\Delta(a)|\,\,
  Q_s(a)\,\,\chi^B_{\gamma}\bigl(a((k_1a)^{-1}\tilde{k})\bigr)\\
  \times\,(v_{\epsilon},\,\pi_W(k_1a)\,
  \pi_W\bigl(k((k_1a)^{-1}\tilde{k})\bigr)\,v_\epsilon)\,\,.
 \ee
Changing the integration variable $k_1\to\tilde{k}k_1$ results in
 \be\label{Qphi3}
  \tilde{\phi}_\epsilon(\tilde{g})\,
  =\frac{1}{|\CW^M|}\,\,\chi^B_{\gamma}\bigl(a(\tilde{g})\bigr)\!\!
  \int\limits_{O_{\ell+1}}\!dk_1\int\limits_{A}\!da\,|\Delta(a)|\,
  Q_s(a)\,\,\chi^B_{\gamma}\bigl(a((k_1a)^{-1})\bigr)\\
  \times\,(v_{\epsilon},\,\pi_W(\tilde{k})\,
  \pi_W\bigl(k_1a\,k((k_1a)^{-1})\bigr)\,v_\epsilon)\\
  =\,\frac{1}{|\CW^M|}\,\chi^B_{\gamma}\bigl(a(\tilde{g})\bigr)\!\!
  \int\limits_{O_{\ell+1}}\!dk_1\int\limits_A\!da\,|\Delta(a)|\,
  Q_s(a)\,\,\chi^B_{\gamma}\bigl(a((k_1a)^{-1})\bigr)\\
  \times\sum_{\e'\in\CP\atop|\e'|=|\e|}(v_{\epsilon},\,\pi_W(\tilde{k})\,v_{\e'})\,
  (v_{\e'},\,\pi_W\bigl(k_1a\,k((k_1a)^{-1})\bigr)\,v_\epsilon)\,\,,
 \ee
where in the latter equality the first identity from \eqref{Wmatel} is
applied.

Now, let us introduce the following functions, elements of
the representation space $\CV_{\epsilon',\gamma}$:
 \be\label{Qmatel}
  \phi^{\epsilon}_{\epsilon',\gamma}(\tilde{g})\,
  =\,(v_{\epsilon},\,\pi_W\bigl(k(\tilde{g})\bigr)\,v_{\epsilon'})\,
  \chi^B_{\gamma}(a(\tilde{g}))\,,
 \ee
and re-write \eqref{Qphi3} as follows:
 \be\label{Qphi4}
  \tilde{\phi}_\epsilon(\tilde{g})\,
  =\,\sum_{\e'\in\CP\atop|\e'|=|\e|}
  \Lambda_{\epsilon,\epsilon'}(s|\gamma)\,\,
  \phi^{\epsilon}_{\epsilon',\gamma}(\tilde{g})\,,
 \ee
where
 \be\label{Qla}
  \Lambda_{\epsilon,\epsilon'}(s|\gamma)\,
  =\frac{1}{|\CW^M|}
  \!\!\int\limits_{O_{\ell+1}\times A}\!\!dk_1\,da\,|\Delta(a)|\,
  Q_s(a)\,\chi^B_{\gamma}\bigl(a((k_1a)^{-1})\bigr)\\
  \times\,(v_{\e'},\,\pi_W\bigl(k_1a\,k((k_1a)^{-1})\bigr)\,v_{\e})\,.
  \ee
Consider the following tautological transformation of \eqref{Qla}
 \be\label{Qla1}
  \Lambda_{\epsilon,\epsilon'}(s|\gamma)\,
  =\frac{1}{|\CW^M|}\,\frac{1}{|M|}\,
  \!\!\sum_{m\in M}\,\int\limits_{O_{\ell+1}\times A}\!\!dk_1\,da\,|\Delta(a)|\,
  Q_s(a)\,\chi^B_{\gamma}\bigl(a((k_1a)^{-1})\bigr)\\
  \times\,(v_{\e'},\,\pi_W\bigl(k_1a\,k((k_1a)^{-1})\bigr)\,v_{\e})\,,
  \ee
and let us make a change of variables $k_1\to mk_1$ in \eqref{Qla1}
 \be\label{Qla11}
  \Lambda_{\epsilon,\epsilon'}(s|\gamma)\,
  =\frac{1}{|\CW^M|}\,\frac{1}{|M|}\,
  \!\!\sum_{m\in M}\,\int\limits_{O_{\ell+1}\times A}\!\!dk_1\,da\,|\Delta(a)|\,
  Q_s(a)\,\chi^B_{\gamma}\bigl(a((mk_1a)^{-1})\bigr)\\
  \times\,(v_{\e'},\,\pi_W\bigl(mk_1a\,k((mk_1a)^{-1})\bigr)\,v_{\e})\,.
  \ee
Using the identities
  \be
a((mk_1a)^{-1})=a(k_1a)^{-1}), \qquad
k((mk_1a)^{-1})=k(k_1a)^{-1})\,m^{-1}\,,
  \ee
we obtain
 \be\label{Qla2}
  \Lambda_{\epsilon,\epsilon'}(s|\gamma)\,
  =\frac{1}{|\CW^M|}\,\int\limits_{O_{\ell+1}\times A}\!\!dk_1\,da\,|\Delta(a)|\,
  Q_s(a)\,\chi^B_{\gamma}\bigl(a((k_1a)^{-1})\bigr)\\
  \times\,
  \frac{1}{|M|}\,  \!\!\sum_{m\in M}
  (v_{\e'},\,\pi_W\bigl(mk_1a\,k((k_1a)^{-1})m^{-1}\bigr)\,v_{\e})\,.
 \ee
Applying the second identity in \eqref{Wmatel} and using the orthogonality
relations for the characters of the finite group $M$ entails
 \be\label{QLambda}
  \Lambda_{\epsilon,\epsilon'}(s|\gamma)\,
  =\delta_{\e,\e'}\,\,
  \frac{1}{|\CW^M|}\,\int\limits_{O_{\ell+1}\times A}\!\!dk_1\,da\,|\Delta(a)|\,
  Q_s(a)\,\chi^B_{\gamma}\bigl(a((k_1a)^{-1})\bigr)\\
  \times\,
  (v_{\e},\,\pi_W\bigl(k_1a\,k((k_1a)^{-1})\bigr)\,v_{\e})\,.
 \ee
Therefore,  \eqref{Qphi4} reduces to the following relation
 \be\label{Qphi5}
  \tilde{\phi}_\epsilon(\tilde{g})\,
  =\,\Lambda_{\e, \e}(s|\gamma)\,\phi_{\epsilon}(\tilde{g})\,\,,
  \ee
with
 \be\label{QLambda1}
  \Lambda_{\e, \e}(s|\gamma)=\,\,
  \frac{1}{|\CW^M|}\,\int\limits_{O_{\ell+1}\times A}\!\!dk_1\,da\,|\Delta(a)|\,
  Q_s(a)\,\chi^B_{\gamma}\bigl(a((k_1a)^{-1})\bigr)\\
  \times\,
 (v_{\e},\,\pi_W\bigl(k_1a\,k((k_1a)^{-1})\bigr)\,v_{\e})\,,
 \ee
i.e.  the $\e$-spherical vector $\phi_{\e}$
is a $\wh{Q}_s$-eigenfunction. To  calculate the
eigenvalue $\Lambda_{\e,\e}$ in \eqref{Qphi5} we use the following
identity
 \be
 k_1a\cdot k((k_1a)^{-1})\,
  =\,n((k_1a)^{-1})^{-1}\cdot a((k_1a)^{-1})^{-1}\,,
 \ee
which gives rise to the  following expression for the matrix element in
\eqref{QLambda1}:
 \be
  (v_{\e},\,\pi_W\bigl(k_1a\,k((k_1a)^{-1})\bigr)\,v_{\e})\,
  =\,\bigl(v_{\epsilon},\,\pi_W\bigl(n((k_1a)^{-1})^{-1}\,
  a((k_1a)^{-1})^{-1}\bigr)\,v_\epsilon\bigr)\\
  =\,\chi^B_{-\i\epsilon}(a((k_1a)^{-1}))\,
  \bigl(v_{\epsilon},\,\pi_W\bigl(n((k_1a)^{-1})^{-1}\bigr)\,v_\epsilon\bigr)\,.
 \ee
Hence, we obtain the following expression for the eigenvalue:
 \be
  \Lambda_{\epsilon,\epsilon}(s|\gamma)\,
  =  \frac{1}{|\CW^M|}\,\int\limits_{O_{\ell+1}\times A}\!\!dk_1\,da\,|\Delta(a)|\,
 Q_s(a)\,\,\chi^B_{\gamma-\i\e}\bigl(a((k_1a)^{-1})\bigr)\\
 \times\,(v_{\e},\,\pi_W\bigl(n((k_1a)^{-1})^{-1}\bigr)\,v_{\e})\,.
 \ee
This can be equivalently written as follows
 \be
  \Lambda_{\epsilon,\epsilon}(s|\gamma)\,
  =\,\frac{1}{|\CW^M|}\!
  \int\limits_{O_{\ell+1}\times A\times O_{\ell+1}}
  \!\!dk_1\,da\,dk_2\,\,|\Delta(a)|\,
  Q_s(a)\,\,\chi^B_{\gamma-\i\e}\bigl(a((k_1ak_2)^{-1})\bigr)\\
  \times\,(v_{\e},\,\pi_W\bigl(n((k_1ak_2)^{-1})^{-1}\bigr)\,v_{\e})\,,
 \ee
and therefore we arrive at the following integral representation for
 the eigenvalue
\be
  \Lambda_{\epsilon,\epsilon}(s|\gamma)\,
  =\!\!\int\limits_{GL_{\ell+1}(\IR)}\!\!d\mu(g)\,\,
  Q_s(g)\,\,\chi_{\gamma-\i\epsilon}^B(g^{-1})
  (v_{\epsilon},\,\pi_W(n(g^{-1}))\,v_\epsilon)\,.
 \ee
Finally using the Iwasawa decomposition $g^{-1}=kan$  we obtain
 \be
  \Lambda_{\epsilon,\epsilon}(s|\gamma)\,
  =\!\!\int\limits_{O_{\ell+1}\times A\times N}\!\!dk\,da\,dn\,\,
  Q_s(g)\,\,\chi^B_{\gamma-\i\epsilon}(a)\,
  (v_{\epsilon},\,\pi_W(n)\,v_\epsilon)\\
  =\!\!\int\limits_{O_{\ell+1}\times A\times N}\,
  dk\,da\,dn\,\,
  Q_s(g)\,\,\chi^B_{\gamma-\i\epsilon}(a)\,.
 \ee
The latter integral is already calculated  in the
 proof of Proposition \ref{SphHB} and the result is in agreement with
\eqref{eL1}.

\noindent {\small {\bf A.A.G.} {\sl Laboratory for Quantum Field
Theory
and Information},\\
\hphantom{xxxx} {\sl Institute for Information
Transmission Problems, RAS, 127994, Moscow, Russia};\\
\hphantom{xxxx} {\it E-mail address}: {\tt anton.a.gerasimov@gmail.com}}\\
\noindent{\small {\bf D.R.L.} {\sl Laboratory for Quantum Field
Theory
and Information},\\
\hphantom{xxxx}  {\sl Institute for Information
Transmission Problems, RAS, 127994, Moscow, Russia};\\
\hphantom{xxxx} {\it E-mail address}: {\tt lebedev.dm@gmail.com}}\\
\noindent{\small {\bf S.V.O.} {\sl
 Beijing Institute of Mathematical Sciences and Applications\,,\\
\hphantom{xxxx} Huairou District, Beijing 101408, China};\\
\hphantom{xxxx} {\it E-mail address}: {\tt oblezin@gmail.com}}

\end{document}